\documentclass[11pt,a4paper]{article}
\usepackage{amsfonts}
\usepackage{latexsym}
\usepackage{cite}
\usepackage{amsmath,amsfonts,latexsym,amssymb}
\usepackage[mathscr]{eucal}
\usepackage{cases,color}
\usepackage{amsthm}
\usepackage{hyperref}
\usepackage[bf,small]{caption2}
\usepackage{float}
\usepackage{graphicx}
\usepackage{amsmath}
\usepackage{amssymb}
\usepackage[all]{xy}

\newtheorem{theorem}{theorem}[section]

\newtheorem{cla}[theorem]{Claim}

\newtheorem{cor}[theorem]{Corollary}

\newtheorem{exmp}[theorem]{Example}
\newtheorem{lem}[theorem]{Lemma}

\newtheorem{que}[theorem]{Question}
\newtheorem{rmk}[theorem]{Remark}
\newtheorem{thm}[theorem]{Theorem}

\begin{document}

\title{\vspace{-2cm}\textbf{High-dimensional components of ${\rm SL}(2,\mathbb{C})$-character varieties of prime knots}}
\author{\Large Haimiao Chen}

\date{}
\maketitle

\begin{abstract}
  For a prime knot $K$, we give sufficient conditions for the existence of a component $\mathcal{C}$ of the irreducible
  ${\rm SL}(2,\mathbb{C})$-character variety of $K$ with $\dim\mathcal{C}>1$, and give a lower bound for $\dim\mathcal{C}$.
  Specifically, we improve a result of Paoluzzi and Porti on Montesinos knots, and positively answer a question posed by Culler and Dunfield in 2018.

\end{abstract}

\section{Introduction}

Let $G={\rm SL}(2,\mathbb{C})$. For a finitely presented group $\Gamma$, call a homomorphism $\rho:\Gamma\to G$ a {\it $G$-representation}.
Define the {\it character} of $\rho$ to be the function $\chi_\rho:\Gamma\to\mathbb{C}$, $x\mapsto{\rm tr}(\rho(x))$.
Call $\rho$ {\it reducible} if the elements of ${\rm Im}(\rho)$ have a common eigenvector; otherwise, call $\rho$ {\it irreducible}.
Call representations $\rho,\rho'$ {\it conjugate} and denote $\rho\sim\rho'$ if there exists $\mathbf{a}\in G$ such that $\rho'(x)=\mathbf{a}\rho(x)\mathbf{a}^{-1}$ for all $x\in\Gamma$.
It is known \cite{CS83} that if $\rho,\rho'$ are irreducible, then $\chi_\rho=\chi_{\rho'}$ if and only if $\rho\sim\rho'$.
Call $\mathcal{R}(\Gamma):=\hom(\Gamma,G)$ the {\it $G$-representation variety} of $\Gamma$.
Let $\mathcal{R}^{\rm irr}(\Gamma)=\{\text{irreducible\ }\rho\in\mathcal{R}(\Gamma)\}$.
The set $\mathcal{X}(\Gamma)=\{\chi_\rho\colon\rho\in\mathcal{R}(\Gamma)\}$ turns out to be an algebraic set, and is called the
$G$-{\it character variety} of $\Gamma$. The subset $\mathcal{X}^{\rm irr}(\Gamma)=\{\chi_\rho\colon\rho\in\mathcal{R}^{\rm irr}(\Gamma)\}$ is Zariski open in $\mathcal{X}(\Gamma)$, and is called the {\it irreducible character variety}.

For a $3$-manifold $M$, abbreviate $\mathcal{X}(\pi_1(M))$ to $\mathcal{X}(M)$.
For a knot $K\subset S^3$, let $E_K$ denote its exterior, and let $\pi(K)=\pi_1(E_K)$; abbreviate $\mathcal{R}(\pi(K))$ to $\mathcal{R}(K)$ and call it the $G$-representation variety of $K$, and so forth.

Within the character variety of a knot, we have special interests in high-dimensional components (HDCs for short), meaning those of dimension $>1$, for several reasons. Firstly, by the Culler-Shalen theory \cite{CS83}, for a 3-manifold $M$, each ideal point of an algebraic curve in $\mathcal{X}(M)$ is associated to an incompressible surface in $M$. If $\dim\mathcal{X}(K)>1$, then a closed incompressible surface in $E_K$ is detected (see \cite[Proposition 2.4]{CCGLS94}).
Note that if $M$ is a closed 3-manifold obtained by a Dehn filling of $E_K$ for such $K$, then $\dim\mathcal{X}(M)>0$.
Secondly, the ${\rm SL}(2,\mathbb{C})$ Casson invariant \cite{Cu01} only uses isolated points in the character variety.
It will be amazing if HDCs are taken into account.
Thirdly, a conjecture suggested by mathematical physics asserts a surprising identity of adjoint Reidemeister torsions \cite{PY23}, and calls for a better understanding of HDC. Fourthly, HDC plays an important role in the study of Kauffman bracket skein module (KBSM) which has become a central object in quantum topology. By \cite{Bu97}, if $\dim\mathcal{X}(M)>0$, then the KBSM of $M$, denoted by $\mathcal{S}(M)$, is not finitely generated. Echoing this, recently it was shown \cite{DKS23} that if $M$ is a closed $3$-manifold with $\dim\mathcal{X}(M)>0$, then $\mathcal{S}(M)$ is not ``tame" (in the sense of \cite[Theorem 3.1 and Proposition 3.2]{DKS23}). In \cite{Ch25} we found a new type of torsion $e$ in $\mathcal{S}(E_K)$ for a $4$-strand Montesinos knot $K$; the proof of $e\ne 0$ was built on investigating a HDC of $\mathcal{X}(K)$. As a related but more general result, \cite[Corollary 2.2]{BD25} states that if $M$ is closed with $\dim\mathcal{X}(M)>0$, then $\mathcal{S}(M)$ contains torsion.

However, by now HDC still remains mysterious. Other than \cite{Kl91,CL96,PP13,PY23}, few results can be found in the literature.

For almost trivial reasons, the character varieties of connected sums of knots always contain HDCs. In this paper, we focus on primes knots.

We give sufficient conditions for a knot $K$ to admit HDCs, meaning that $\mathcal{X}^{\rm irr}(K)$ contains HDCs, and give lower bounds for the dimensions (see Theorem \ref{thm:HDC}). Then specifically, we improve a result of Paoluzzi and Porti on the dimensions of character varieties of a class of Montesinos knots (see Remark \ref{rmk:improve}).
We also find two kinds of ``conditional" substructures for $K$ to admit HDCs, in Section \ref{sec:conditional}.
Along the way, we clarify some useful facts that were never mentioned by other people, using elementary techniques.

We work out details for the $(3,3,3,4)$-pretzel knot and two $12$-crossing knots, all of which are hyperbolic.
In particular, we positively answer the following question posed by Culler and Dunfield in 2018:
\begin{que}[Question (7) in {\cite[Section 9]{CD18}}]\label{que:CD-2018}
Does there exist a closed atoroidal $3$-manifold $M$ such that $\dim H_1(M;\mathbb{Q})\le 1$ and
$\dim\mathcal{X}^{\rm irr}(M)\ge 1$? What if one restricts to 0-surgery on a knot in $S^3$?
\end{que}


\noindent
{\bf Acknowledgements}

I'd like to thank Professor Jiming Ma for inspiring conversations when I visited Fudan University, and thank Doctor Seokbeom Yoon for valuable discussions. I also thank my colleague Professor Yubo Luo for improving the working environment.

\section{Preliminary}

Most notations adopted in this paper had been introduced in \cite[Section 2]{Ch22}. 
One may refer there if necessary.

Use bold letters to denote elements of $G$.
For $\mathbf{a},\mathbf{b}\in G$, let $\mathbf{a}\lrcorner\mathbf{b}=\mathbf{a}\mathbf{b}\mathbf{a}^{-1}$.

For $\kappa\in\mathbb{C}^\ast=\mathbb{C}\setminus\{0\}$, $u\in\mathbb{C}$, put
$$\mathbf{d}(\kappa)=\left(\begin{array}{cc} \kappa & 0 \\ 0 & \kappa^{-1} \end{array}\right), \qquad
\mathbf{u}(\kappa)=\left(\begin{array}{cc} \kappa & 1 \\ 0 & \kappa^{-1} \end{array}\right), \qquad
\mathbf{p}(u)=\left(\begin{array}{cc} 1 & u \\ 0 & 1 \end{array}\right).$$
Let $\mathbf{e}=\mathbf{d}(1)$, the identity matrix, and let $\mathbf{p}=\mathbf{p}(1)$.

For $t\in\mathbb{C}$, let $\mathbb{B}_t=\mathbb{C}\setminus\{2,t^2-2\}$, let
$G(t)=\{\mathbf{x}\in G\setminus\{\pm\mathbf{e}\}\colon {\rm tr}(\mathbf{x})=t\}.$

Each element of $G(t)$ with $t\ne\pm2$ is conjugate to $\mathbf{d}(\kappa)$ for some $\kappa$ with $\kappa+\kappa^{-1}=t$, and each element of
$G(2\varepsilon)$ with $\varepsilon\in\{\pm1\}$ is conjugate to $\varepsilon\mathbf{p}$.

For $\mathbf{a}\in G$, let ${\rm Cen}(\mathbf{a})=\{\mathbf{c}\in G\colon \mathbf{c}\mathbf{a}=\mathbf{a}\mathbf{c}\}$, the centralizer of $\mathbf{a}$.

Given $\mathbf{a}\in G(t)$, $r\in\mathbb{B}_t$, let $\mathcal{C}_t^r(\mathbf{a})=\{\mathbf{x}\in G(t)\colon {\rm tr}(\mathbf{a}\mathbf{x})=r\}$.

\begin{lem}\label{lem:centralizer}
{\rm(i)} For each $\mathbf{a}\in G(t)$, we have ${\rm Cen}(\mathbf{a})\cong\mathbb{C}^\ast$ if $t\ne \pm 2$, and
${\rm Cen}(\mathbf{a})\cong\mathbb{C}\sqcup\mathbb{C}$ if $t\in\{\pm 2\}$. In any case,
$${\rm Cen}(\mathbf{a})=\{\mu\mathbf{a}+\nu\mathbf{e}\colon\mu^2+t\mu\nu+\nu^2=1\}.$$

{\rm(ii)} If $\mathbf{a}\mathbf{b}\ne\mathbf{b}\mathbf{a}$, then ${\rm Cen}(\mathbf{a})\cap{\rm Cen}(\mathbf{b})=\{\pm\mathbf{e}\}$.
\end{lem}

\begin{proof}
(i) Up to conjugacy we may assume $\mathbf{a}=\mathbf{d}(\kappa)$ with $\kappa\ne\pm1$ or $\mathbf{a}=\varepsilon\mathbf{p}$ with $\varepsilon\in\{\pm1\}$. In the former case, ${\rm Cen}(\mathbf{a})=\{\mathbf{d}(\mu)\colon\mu\in\mathbb{C}^\ast\}$; in the latter case,
${\rm Cen}(\mathbf{a})=\{\pm\mathbf{p}(u)\colon u\in\mathbb{C}\}$. Then the assertions are easy to verify.

(ii) By (i), each $\mathbf{c}\in{\rm Cen}(\mathbf{a})$ has the form $\mu\mathbf{a}+\nu\mathbf{e}$. If furthermore $\mathbf{c}\in{\rm Cen}(\mathbf{b})$, then
$\mu\mathbf{a}\mathbf{b}=\mu\mathbf{b}\mathbf{a}$, hence $\mu=0$, implying $\mathbf{c}\in\{\pm\mathbf{e}\}$.
\end{proof}

\begin{lem}\label{lem:basic}
For any $\mathbf{a}\in G(t)$, $r\in\mathbb{B}_t$, we have $\mathcal{C}_t^r(\mathbf{a})\cong \mathbb{C}^\ast$ if $t\ne\pm2$, and $\mathcal{C}_t^r(\mathbf{a})\cong\mathbb{C}$ if $t\in\{\pm2\}$.
\end{lem}

\begin{proof}
Each $\mathbf{x}\in\mathcal{C}_t^r(\mathbf{a})$ can be written as
$$\mathbf{x}=\left(\begin{array}{cc} t/2+a & b \\ c & t/2-a \end{array}\right), \qquad bc=\frac{t^2}{4}-1-a^2.$$

When $t\ne\pm2$, up to conjugacy we may assume $\mathbf{a}=\mathbf{d}(\kappa)$ with $\kappa+\kappa^{-1}=t$.
It follows from ${\rm tr}(\mathbf{a}\mathbf{x})=r$ that $(\kappa-\kappa^{-1})a=r-t^2/2$, and then
$$bc=\frac{t^2}{4}-1-\Big(\frac{r-t^2/2}{\kappa-\kappa^{-1}}\Big)^2=\frac{(r-2)(r-t^2+2)}{4-t^2}\ne 0.$$

When $t=2\varepsilon$ with $\varepsilon\in\{\pm1\}$, up to conjugacy we may assume $\mathbf{a}=\varepsilon\mathbf{p}$.
From ${\rm tr}(\mathbf{a}\mathbf{x})=r$ we see $c=\varepsilon(r-2)\ne 0$, then $b$ is determined by $a$ via $b=-a^2/c$.
\end{proof}

Call a tuple $\underline{\mathbf{a}}=(\mathbf{a}_1,\ldots,\mathbf{a}_k)\in G^k$ {\it reducible} if $\mathbf{a}_1,\ldots,\mathbf{a}_k$ share an eigenvector; call $\underline{\mathbf{a}}$ {\it irreducible} otherwise.
Say that $\underline{\mathbf{a}}'=(\mathbf{a}'_1,\ldots,\mathbf{a}'_k)$ is {\it conjugate} to $\underline{\mathbf{a}}$ if there exists $\mathbf{c}\in G$ with $\mathbf{c}\lrcorner\underline{\mathbf{a}}=\underline{\mathbf{a}}'$, i.e. $\mathbf{c}\lrcorner\mathbf{a}_i=\mathbf{a}'_i$ for $1\le i\le k$.

Let
\begin{align*}
f_{t}(r_1,r_2,r_3;r)=\ &r^2+t(t^2-r_1-r_2-r_3)r+t^2(3-r_1-r_2-r_3)   \\
&+r_1^2+r_2^2+r_3^2+r_1r_2r_3-4.
\end{align*}

The following was stated in \cite[Lemma 3.1]{Ch25}; there, actually $t=\pm 2$ can be included. One may refer to \cite[Section 2 and 5]{Go09}.
\begin{lem}\label{lem:matrix}
{\rm(i)} $(\mathbf{a}_1,\mathbf{a}_2)\in G(t)^2$ is irreducible if and only if ${\rm tr}(\mathbf{a}_1\mathbf{a}_2)\in\mathbb{B}_t$.

{\rm(ii-1)} Given $t_{12}\in\mathbb{B}_t$, up to conjugacy there exists a unique pair $(\mathbf{a}_1,\mathbf{a}_2)\in G(t)^2$ with
${\rm tr}(\mathbf{a}_1\mathbf{a}_2)=t_{12}$.

{\rm(ii-2)} If furthermore $t_{13},t_{23},t_{123}$ are given with $f_{t}(t_{12},t_{13},t_{23};t_{123})=0$, then there exists a unique $\mathbf{a}_3\in G(t)$ such that ${\rm tr}(\mathbf{a}_1\mathbf{a}_3)=t_{13}$, ${\rm tr}(\mathbf{a}_2\mathbf{a}_3)=t_{23}$,
${\rm tr}(\mathbf{a}_1\mathbf{a}_2\mathbf{a}_3)=t_{123}$.
\end{lem}

\medskip

For a knot $K$, in the sprit of Wirtinger presentation, a representation of $K$ identifies with a map $\rho$ assigning to each directed arc $\mathsf{a}$ an element $\rho(\mathsf{a})\in G$ such that $\rho(\mathsf{a}^{-1})=\rho(\mathsf{a})^{-1}$ (where $\mathsf{a}^{-1}$ is obtained by reversing the direction of $\mathsf{a}$), and $\rho(\mathsf{c})=\rho(\mathsf{a})\lrcorner\rho(\mathsf{b})$ for any $\mathsf{a},\mathsf{b},\mathsf{c}$ forming a crossing as
\begin{figure} [h]
  \centering
  \includegraphics[width=2cm]{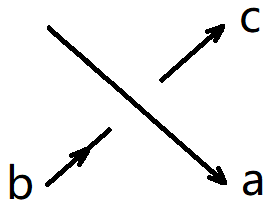}. \\
\end{figure}
\\
To present such $\rho$, it suffices to choose a direction for each arc and label this arc with an element of $G$.
For $t\in\mathbb{C}$, let $\mathcal{R}_t(K)$ denote the set of $\rho$'s such that $\rho(\mathsf{a})\in G(t)$ for each directed arc $\mathsf{a}$, and
let $\mathcal{R}_t^{\rm irr}(K)$ denote the subset of the irreducible ones; 
let $\mathcal{X}_t^{\rm irr}(K)=\{\chi_\rho\colon\rho\in\mathcal{R}_t^{\rm irr}(K)\}$, which identifies with the set of conjugacy classes of
$\rho\in\mathcal{R}_t^{\rm irr}(K)$.

In this way, we can talk about representations of tangles. Given a tangle $T$, let $\mathcal{R}(T)$ denote the set of all representations of $T$, and so forth.

Given $\rho\in\mathcal{R}(T)$ and $\mathbf{c}\in G$, let $\mathbf{c}\lrcorner\rho\in\mathcal{R}(T)$ send each directed arc $\mathsf{a}$ to $\mathbf{c}\lrcorner\rho(\mathsf{a})$.
Call $\rho,\rho'$ {\it conjugate} and denote $\rho\sim\rho'$ if $\rho'=\mathbf{c}\lrcorner\rho$ for some $\mathbf{c}$.

If $S$ is a subtangle of $T$, then by restriction each $\rho\in\mathcal{R}(T)$ gives rise to a representation of $S$, which is denoted by $\rho|_{S}$.

\begin{figure}[H]
  \centering
  \includegraphics[width=12cm]{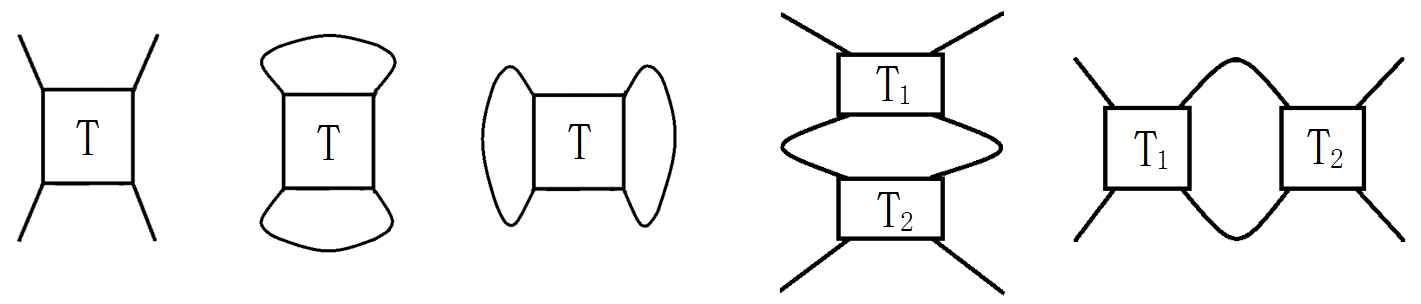}\\
  \caption{From left to right: $T\in\mathcal{T}^2_2$; $N(T)$; $D(T)$; $T_1\ast T_2$; $T_1+T_2$.}\label{fig:tangle}
\end{figure}

Let $\mathcal{T}^2_2$ denote the set of tangles with four ends. For $T\in\mathcal{T}^2_2$, let $T^{\rm nw}$, $T^{\rm ne}$, $T^{\rm sw}$, $T^{\rm se}$ respectively denote the northwest-, northeast-, southwest- and southeast ends, all directed outward.
From $T$ one may construct two links: the numerator $N(T)$ and the denominator $D(T)$. Defined on $\mathcal{T}^2_2$ are the vertical composition $\ast$ and horizontal composition $+$.
See Figure \ref{fig:tangle}.

A positive (reps. negative) crossing is regarded as a tangle, denoted by $[1]$ (resp. $[-1]$).
For $k\ne 0$, the horizontal composite of $|k|$ copies of $[1]$ (resp. $[-1]$) is denoted by $[k]$ if $k>0$ (resp. $k<0$), and the
vertical composite of $|k|$ copies of $[1]$ (resp. $[-1]$) is denoted by $[1/k]$ if $k>0$ (resp. $k<0$).

To each $p/q\in\mathbb{Q}$ is associated a rational tangle
$$[p/q]:=\begin{cases}
[k_{1}]\ast [1/k_{2}]+\cdots+[k_{s}], &2\nmid s \\
[1/k_{1}]+[k_{2}]\ast\cdots +[k_{s}], &2\mid s
\end{cases},$$
if $[k_s;k_{s-1},\ldots,k_1]$ is a continued fraction for $p/q$. 
Let $B(p/q=N([p/q]))$, a rational link.
Given $p_1/q_1,\ldots,p_m/q_m\in\mathbb{Q}$, denote $D([p_1/q_1]\ast\cdots\ast[p_m/q_m])$ by $M(p_1/q_1,\ldots,p_m/q_m)$ and call it a Montesinos link.

Given $\rho\in\mathcal{R}(T)$ for $T\in\mathcal{T}^2_2$, let
\begin{alignat*}{2}
\mathbf{g}_\rho&=\rho(T^{\rm nw})\rho(T^{\rm ne}), \qquad  &{\rm tr}_h(\rho)&={\rm tr}(\mathbf{g}_\rho);   \\
\dot{\mathbf{g}}_\rho&=\rho(T^{\rm sw})\rho(T^{\rm nw}),  \qquad  &{\rm tr}_v(\rho)&={\rm tr}(\dot{\mathbf{g}}_\rho).
\end{alignat*}

\section{Structures for high-dimensional components}

\subsection{A general structure}

\begin{lem}\label{lem:moduli}
Given $\mathbf{a},\mathbf{b}\in G(t)$ and $\underline{a}=(a_1,\ldots,a_n)\in\mathbb{B}_t^{n}$ with $n\ge 3$, put
$$\mathcal{M}_t^{\underline{a}}(\mathbf{a},\mathbf{b})=\big\{(\mathbf{x}_1,\ldots,\mathbf{x}_{n-1})\in G(t)^{n-1}\colon
{\rm tr}(\mathbf{x}_{i-1}\mathbf{x}_i)=a_i, 1\le i\le n\big\},$$
where $\mathbf{x}_0=\mathbf{a}$, $\mathbf{x}_n=\mathbf{b}$ for convenience.
Then $\dim\mathcal{M}_t^{\underline{a}}(\mathbf{a},\mathbf{b})=n-2$.
\end{lem}

\begin{proof}
Use induction on $n$ to prove the assertion.

First, suppose $n=3$. Define
$$\Phi:\mathcal{M}_t^{\underline{a}}(\mathbf{a},\mathbf{b})\to\mathcal{C}_t^{a_1}(\mathbf{a}), \qquad  (\mathbf{x}_1,\mathbf{x}_2)\mapsto\mathbf{x}_1.$$
\begin{enumerate}
  \item If $(\mathbf{x}_1,\mathbf{b})$ is irreducible, then by Lemma \ref{lem:matrix} (ii-2), 
        there are two or one matrix $\mathbf{x}_2\in G(t)$ with ${\rm tr}(\mathbf{x}_1\mathbf{x}_2)=a_2$ and
        ${\rm tr}(\mathbf{x}_2\mathbf{b})=a_3$.
  \item If $(\mathbf{x}_1,\mathbf{b})$ is reducible and $\mathbf{x}_1\mathbf{b}\ne\mathbf{b}\mathbf{x}_1$,
        then $t\ne\pm2$, and up to conjugacy we may assume
        $\mathbf{x}_1=\mathbf{d}(\kappa)$, $\mathbf{b}=\mathbf{u}(\kappa^{\varepsilon})$,
        with $\kappa+\kappa^{-1}=t$ and $\varepsilon\in\{\pm1\}$.
        For $\mathbf{x}_2=(y_{ij})_{2\times 2}\in G(t)$, the conditions
        ${\rm tr}(\mathbf{x}_1\mathbf{x}_2)=a_2$ and ${\rm tr}(\mathbf{b}\mathbf{x}_2)=a_3$ are equivalent to
        \begin{align*}
        y_{11}=\frac{a_2-\kappa^{-1}t}{\kappa-\kappa^{-1}}, \qquad y_{22}&=\frac{\kappa t-a_2}{\kappa-\kappa^{-1}}, \qquad
        y_{21}=a_3-\varepsilon a_2+\frac{\varepsilon-1}{2}t^2, \\
        y_{12}y_{21}&=y_{11}y_{22}-1=\frac{(a_2-2)(a_2-t^2+2)}{4-t^2}.
        \end{align*}
        Hence when $a_3\ne \varepsilon a_2+(1-\varepsilon)t^2/2$, there exists a unique $\mathbf{x}_2\in G(t)$ with
        ${\rm tr}(\mathbf{x}_1\mathbf{x}_2)=a_2$, ${\rm tr}(\mathbf{b}\mathbf{x}_2)=a_3$; otherwise,
        no such $\mathbf{x}_2$ exists.
  \item If $\mathbf{x}_1\mathbf{b}=\mathbf{b}\mathbf{x}_1$, so that $\mathbf{x}_1=\mathbf{b}^\varepsilon$ with $\varepsilon\in\{\pm 1\}$, 
        then it is easy to see that the space of $\mathbf{x}_2$ with ${\rm tr}(\mathbf{x}_1\mathbf{x}_2)=a_2$ and ${\rm tr}(\mathbf{b}\mathbf{x}_2)=a_3$ is
        $\mathcal{C}_t^{a_3}(\mathbf{b})$ if $a_2=\varepsilon a_3+(1-\varepsilon)t^2/2$, and $\emptyset$ otherwise.
\end{enumerate}
By Lemma \ref{lem:basic}, $\dim\mathcal{C}_t^{a_1}(\mathbf{a})=\dim\mathcal{C}_t^{a_3}(\mathbf{b})=1$. Therefore, we have shown that $\dim\Phi^{-1}(\mathcal{C}_t^{a_1}(\mathbf{a})\setminus\{\mathbf{b}^{\pm1}\})=1$ and $\dim\Phi^{-1}(\mathbf{b}^{\pm1})\le 1$,
so $\dim\mathcal{M}_t^{\underline{a}}(\mathbf{a},\mathbf{b})=1$.

Now suppose $n=k+1\ge 4$ and that the assertion holds for $n=k$. For $\underline{a}=(a_1,\ldots,a_{k+1})$, let $\underline{a}'=(a_1,\ldots,a_{k})$, and define
$$\mathcal{M}_t^{\underline{a}}(\mathbf{a},\mathbf{b})\to\mathcal{C}_t^{a_{k+1}}(\mathbf{b}), \qquad
(\mathbf{x}_1,\ldots,\mathbf{x}_k)\mapsto \mathbf{x}_k.$$
The pre-image of each $\mathbf{x}\in\mathcal{C}_t^{a_{k+1}}(\mathbf{b})$ is $\mathcal{M}_t^{\underline{a}'}(\mathbf{a},\mathbf{x})$, which is $(k-2)$-dimensional by the inductive hypothesis.
Thus, $\dim\mathcal{M}_t^{\underline{a}}(\mathbf{a},\mathbf{b})=k-1=n-2$.
\end{proof}

\begin{thm}\label{thm:HDC}
Suppose $K=D(T_0\ast\cdots\ast T_{n})$, with $T_0,\ldots,T_n\in\mathcal{T}^2_2$, $n\ge 3$.
Suppose for each generic $t$ there exist $\varrho_i\in\mathcal{R}_t(T_i)$, $0\le i\le n$ such that
$\mathbf{g}_{\varrho_i}=\mathbf{e}$ for $0\le i\le n$, and ${\rm tr}_v(\varrho_i)\ne 2,t^2-2$ for $1\le i\le n$.
Then
$$\dim\mathcal{X}^{\rm irr}(K)\ge\begin{cases} n-1,&\varrho_0(T_0^{\rm ne})\varrho_0(T_0^{\rm se})\ne\varrho_0(T_0^{\rm ne})\varrho_0(T_0^{\rm se}) \\ n-2,&\varrho_0(T_0^{\rm ne})\varrho_0(T_0^{\rm se})=\varrho_0(T_0^{\rm ne})\varrho_0(T_0^{\rm se}) \end{cases}.$$
\end{thm}

\begin{proof}
Let $t\in\mathbb{C}$ be generic; let $\underline{a}=(a_1,\ldots,a_n)$, with $a_i=t^2-{\rm tr}_v(\varrho_i)$.
Let $\mathbf{a}=\varrho_0(T_0^{\rm se})$, $\mathbf{b}=\varrho_0(T_0^{\rm ne})^{-1}$.

Given any $\underline{\mathbf{x}}=(\mathbf{x}_1,\ldots,\mathbf{x}_{n-1})\in\mathcal{M}_t^{\underline{a}}(\mathbf{a},\mathbf{b})$, let $\mathbf{x}_0=\mathbf{a}$, $\mathbf{x}_n=\mathbf{b}$.

For each $1\le i\le n$, since
${\rm tr}(\mathbf{x}_{i-1}^{-1}\mathbf{x}_i)={\rm tr}((t\mathbf{e}-\mathbf{x}_{i-1})\mathbf{x}_i)=t^2-a_i={\rm tr}_v(\varrho_i)$ and $t^2-a_i\in\mathbb{B}_t$,
by Lemma \ref{lem:matrix} (ii-1) there exists $\mathbf{c}_i\in G$ such that
$$\mathbf{c}_i\lrcorner(\varrho_i(T_i^{\rm ne}),\varrho_i(T_i^{\rm se}))=(\mathbf{x}_{i-1}^{-1},\mathbf{x}_i).$$
Put $\rho_i=\mathbf{c}_i\lrcorner\varrho_i\in\mathcal{R}_t(T_i)$. Then
$\rho_i(T_i^{\rm ne})=\mathbf{x}_{i-1}^{-1}$, $\rho_i(T_i^{\rm se})=\mathbf{x}_i$, and $\mathbf{g}_{\rho_i}=\mathbf{e}$.

Let $\rho_0=\varrho_0$. Due to the compatibilities of the values at common arcs:
$$\rho_{i-1}(T_{i-1}^{\rm se})=\rho_i(T_i^{\rm ne})^{-1} \ \ \  (1\le i\le n), \qquad \rho_n(T_n^{\rm se})=\rho_0(T_0^{\rm ne})^{-1},$$
we can glue $\rho_0,\ldots,\rho_n$ into a representation $\rho^{\underline{\mathbf{x}}}\in\mathcal{R}_t^{\rm irr}(K)$.

Let $H={\rm Cen}(\mathbf{a})\cap{\rm Cen}(\mathbf{b})$, and $\overline{H}=H/\{\pm\mathbf{e}\}$. Then $\overline{H}$ freely acts on $\mathcal{M}_t^{\underline{a}}(\mathbf{a},\mathbf{b})$ by conjugation.
The map $\underline{\mathbf{x}}\mapsto\chi_{\rho^{\underline{\mathbf{x}}}}$
factors through the action of $\overline{H}$, to define a map
$$\Psi_t:\mathcal{M}_t^{\underline{a}}(\mathbf{a},\mathbf{b})/\overline{H}\to\mathcal{X}_t^{\rm irr}(K).$$
If $\rho^{\underline{\mathbf{x}}}\sim\rho^{\underline{\mathbf{x}}'}$, say $\rho^{\underline{\mathbf{x}}'}=\mathbf{c}\lrcorner\rho^{\underline{\mathbf{x}}}$ for some $\mathbf{c}\in G$, then $\mathbf{c}\in H$ and $\underline{\mathbf{x}}'=\mathbf{c}\lrcorner\underline{\mathbf{x}}'$. Hence $\Psi_t$ is injective.

If $\mathbf{a}\mathbf{b}\ne\mathbf{b}\mathbf{a}$, then by Lemma \ref{lem:centralizer} (ii), $H=\{\pm\mathbf{e}\}$.
If $\mathbf{a}\mathbf{b}=\mathbf{b}\mathbf{a}$, then $\mathbf{b}\in\{\mathbf{a}^{\pm1}\}$, so $H={\rm Cen}(\mathbf{a})$; by Lemma \ref{lem:centralizer} (i), $\dim H=1$.
Consequently,
$$\dim\mathcal{X}_t^{\rm irr}(K)\ge\dim\mathcal{M}_t^{\underline{a}}(\mathbf{a},\mathbf{b})
-\dim H=\begin{cases} n-2,&\mathbf{a}\mathbf{b}\ne\mathbf{b}\mathbf{a} \\
n-3,&\mathbf{a}\mathbf{b}=\mathbf{b}\mathbf{a} \end{cases}.$$
Then the lower bound for $\dim\mathcal{X}^{\rm irr}(K)$ follows.
\end{proof}

\begin{cor}\label{cor:HDC}
Let $K=M(\beta_0/\alpha_0,\ldots,\beta_n/\alpha_n)$ be a Montesinos knot with $n\ge 3$. If each $\alpha_i\ge 3$, then
$\dim\mathcal{X}^{\rm irr}(K)\ge n-1$.
\end{cor}

\begin{proof}
According to \cite[Section 7]{ORS08} (also see \cite[Section 3.1]{Ch25}), for each generic $t$, there exists $\sigma_i\in\mathcal{R}^{\rm irr}(B(\beta_i/\alpha_i))$, which identifies with some $\varrho_i\in\mathcal{R}_t(T_i)$ such that $\mathbf{g}_{\varrho_i}=\mathbf{e}$ and
${\rm tr}_v(\varrho_i)\ne 2,t^2-2$. In particular, $\varrho_0(T_0^{\rm ne})\varrho_0(T_0^{\rm se})\ne\varrho_0(T_0^{\rm ne})\varrho_0(T_0^{\rm se})$. By Theorem \ref{thm:HDC}, $\dim\mathcal{X}^{\rm irr}(K)\ge n-1$.
\end{proof}

\begin{rmk}\label{rmk:improve}
\rm 
On restriction to the case $2\mid\alpha_0$ (and $2\nmid\alpha_i$ for $1\le i\le n$), This improves the result of
\cite[Theorem 1]{PP13} that $\dim\mathcal{X}(K)\ge n-2$.
The key ingredient is that \cite[Section 4]{PP13} requires $\rho(\mu_1)\rho(\mu_{n+1})=\rho(\mu_{n+1})\rho(\mu_1)$.
The commuting condition decreases the dimension by $1$, as seen in the last paragraph of the proof of Theorem \ref{thm:HDC}.
\end{rmk}

\begin{figure}[h]
  \centering
  \includegraphics[width=9cm]{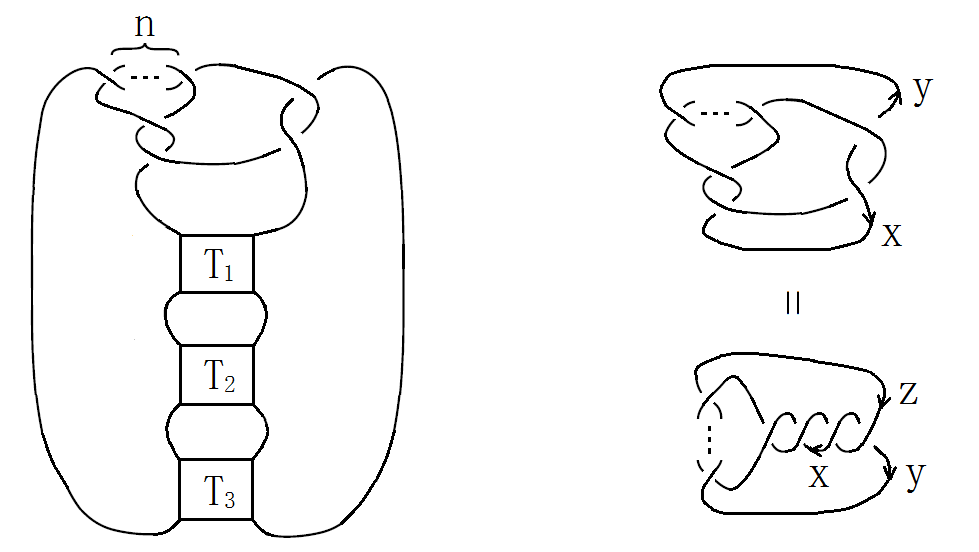}\\
  \caption{Left: the knot $K_n$. Right: $N(T_0)=B(-(4n+1)/n)$.}\label{fig:example-1}
\end{figure}

\begin{exmp}[A non-Montesinos knot]
\rm Let $K_n=D(T_0\ast T_1\ast T_2\ast T_3)$, with $T_0=([n]\ast[1/2])+[1/2]$, and $T_i=[p_i/q_i]$, $i=1,2,3$. See Figure \ref{fig:example-1}, left.
Suppose $q_1,q_2,q_3\ge 3$, and $n\ge 3$.

Note that $N(T_0)=B(-(4n+1)/n)$, as illustrated in Figure \ref{fig:example-1}, right.
For each generic $t$, there exists $\tau\in\mathcal{R}_t^{\rm irr}(N(T_0))$, which is the same as some $\varrho_0\in\mathcal{R}_t(T_0)$ with $\mathbf{g}_{\varrho_0}=\mathbf{e}$. Then $\mathbf{y}:=\tau(y)=\varrho_0(T_0^{\rm ne})$ does not commute with $\mathbf{x}:=\tau(x)=\varrho_0(T_0^{\rm se})$.
This is because $\mathbf{y}=\mathbf{z}\lrcorner\mathbf{x}^{-1}$, where $\mathbf{z}=\tau(z)$; if $\mathbf{x}\mathbf{y}=\mathbf{y}\mathbf{x}$, then $\mathbf{z}\lrcorner\mathbf{x}=\mathbf{x}^{\varepsilon}$ for some $\varepsilon\in\{\pm 1\}$, so that $$\mathbf{z}^2\lrcorner\mathbf{x}=\mathbf{z}\lrcorner(\mathbf{z}\lrcorner\mathbf{x})=\mathbf{z}\lrcorner\mathbf{x}^{\varepsilon}
=(\mathbf{z}\lrcorner\mathbf{x})^\varepsilon=(\mathbf{x}^\varepsilon)^\varepsilon=\mathbf{x},$$
implying $\mathbf{x}\in\{\mathbf{z}^{\pm2}\}$, so $\mathbf{y}\mathbf{z}=\mathbf{z}\mathbf{y}$, which contradicts the irreducibility of $\tau$.

For each $1\le i\le 3$, similarly as in the proof of Corollary \ref{cor:HDC}, there exists $\varrho_i\in\mathcal{R}_t(T_i)$ such that $\mathbf{g}_{\varrho_i}=\mathbf{e}$ and ${\rm tr}_v(\varrho_i)\ne 2,t^2-2$.

By Theorem \ref{thm:HDC}, $\dim\mathcal{X}^{\rm irr}(K_n)\ge 2$.
\end{exmp}

A nontrivial answer to the following will be very interesting:
\begin{que}\label{que:HDC}
Is there a prime knot $K$ not of the form $D(T_0\ast\cdots\ast T_n)$ with $T_i\in\mathcal{T}_2^2$ and $n\ge 3$ that admits a HDC?
\end{que}

\subsection{Conditional structures}\label{sec:conditional}

We partially answer Question \ref{que:HDC}, by showing that when
$$K=D(T_1\ast T_2\ast(T_3+T_4))$$
for some $T_i\in\mathcal{T}_2^2$,
there are two possibilities for $K$ to admit a HDC. In either case, the $T_i$'s are required to satisfy some conditions.

\begin{figure}[H]
  \centering
  \includegraphics[width=9cm]{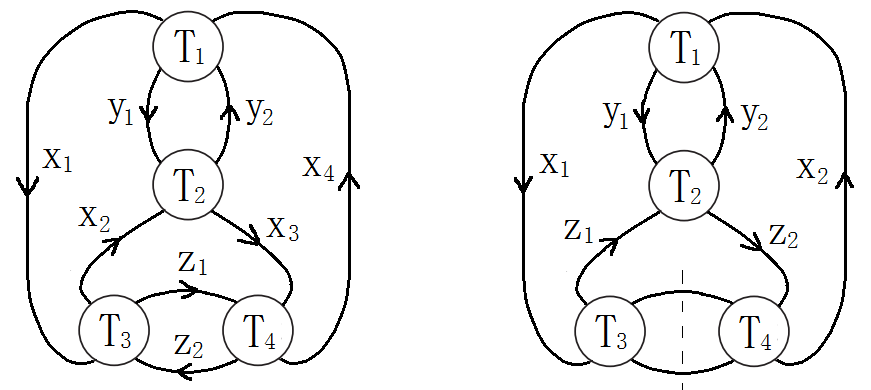} \\
  \caption{}\label{fig:conditional}
\end{figure}

{\bf Case 1}. Suppose that for each generic $t$, there exist $\varrho_i\in\mathcal{R}_t(N(T_i))$, $i=1,2$
and $\varrho_j\in\mathcal{R}_t(D(T_j))$, $j=3,4$ such that
${\rm tr}_v(\varrho_1)={\rm tr}_v(\varrho_2)=:a\in\mathbb{B}_t$, and ${\rm tr}_h(\varrho_3)={\rm tr}_h(\varrho_4)=:b\in\mathbb{B}_t$.
This requires $T_1,T_2$ to be related in a certain way, and similarly for $T_3,T_4$.

Given $(\mathbf{x},\mathbf{y},\mathbf{z})\in G(t)^3$ with ${\rm tr}(\mathbf{x}\mathbf{y})=a$, ${\rm tr}(\mathbf{x}\mathbf{z})=b$, similarly as in the proof of Theorem \ref{thm:HDC}, there exist $\rho_i\in\mathcal{R}_t(T_i)$, $1\le i\le 4$ such that
\begin{alignat*}{4}
\rho_1(T_1^{\rm nw})&=\mathbf{x}, \qquad &\rho_1(T_1^{\rm ne})&=\mathbf{x}^{-1}, \qquad &\rho_1(T_1^{\rm sw})&=\mathbf{y},
\qquad   &\rho_1(T_1^{\rm ne})&=\mathbf{y}^{-1},  \\
\rho_2(T_2^{\rm se})&=\mathbf{x}, \qquad &\rho_2(T_2^{\rm sw})&=\mathbf{x}^{-1}, \qquad &\rho_2(T_2^{\rm ne})&=\mathbf{y},
\qquad   &\rho_2(T_2^{\rm nw})&=\mathbf{y}^{-1},  \\
\rho_3(T_3^{\rm nw})&=\mathbf{x}, \qquad &\rho_3(T_3^{\rm sw})&=\mathbf{x}^{-1}, \qquad &\rho_3(T_3^{\rm ne})&=\mathbf{z},
\qquad   &\rho_3(T_3^{\rm se})&=\mathbf{z}^{-1},  \\
\rho_4(T_4^{\rm se})&=\mathbf{x}, \qquad &\rho_4(T_4^{\rm ne})&=\mathbf{x}^{-1}, \qquad &\rho_4(T_4^{\rm sw})&=\mathbf{z},
\qquad   &\rho_4(T_4^{\rm nw})&=\mathbf{z}^{-1}.
\end{alignat*}
We can glue the $\rho_i$'s into some $\rho\in\mathcal{R}_t^{\rm irr}(K)$ with $\rho(x_i)=\mathbf{x}$ for $1\le i\le 4$ and
$\rho(y_i)=\mathbf{y}$, $\rho(z_i)=\mathbf{z}$ for $i=1,2$.
See the left part of Figure \ref{fig:conditional}. 

By Lemma \ref{lem:matrix}, the conjugacy classes of $(\mathbf{x},\mathbf{y},\mathbf{z})$ are parameterized by $c:={\rm tr}(\mathbf{y}\mathbf{z})$ and
$r:={\rm tr}(\mathbf{x}\mathbf{y}\mathbf{z})$, subject to $f_t(a,b,c;r)=0$.
Therefore, $K$ admits a HDC.

\begin{figure}[H]
  \centering
  \includegraphics[width=6.5cm]{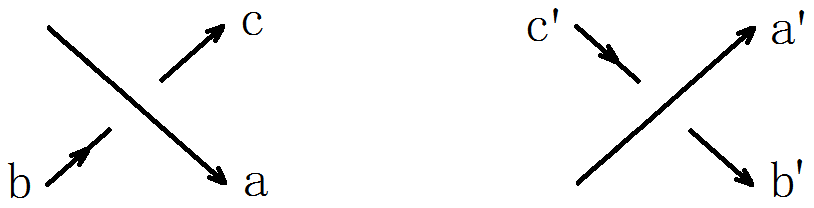}\\
  \caption{Left: a crossing in $T_3$. Right: the corresponding one in $T_4$ under $\sigma$.}\label{fig:symmetric}
\end{figure}

{\bf Case 2}. Suppose $T_4=\sigma(T_3)$, where $\sigma$ is the reflection along the dotted line (see the right part of
Figure \ref{fig:conditional}), so $\mathsf{a}\mapsto\mathsf{a}':=\sigma(\mathsf{a})^{-1}$ establishes a bijection between directed arcs of $T_3$ and those of $T_4$. Furthermore, suppose that for each generic $t$, there exists
$\varrho_i\in\mathcal{R}_t(N(T_i))$ with ${\rm tr}_v(\varrho_i)=a_i\in\mathbb{B}_t$, $i=1,2$, and for each generic $b$, there exists $\tau\in\mathcal{R}_t(T_3)$ with ${\rm tr}_v(\tau)=b$.

Given $(\mathbf{x},\mathbf{y},\mathbf{z})\in G(t)^3$ with ${\rm tr}(\mathbf{x}\mathbf{y})=a_1$, ${\rm tr}(\mathbf{y}\mathbf{z})=a_2$
and ${\rm tr}(\mathbf{x}\mathbf{z})$ being generic, there exist $\rho_i\in\mathcal{R}_t(T_i)$, $i=1,2,3$ such that
\begin{alignat*}{4}
\rho_1(T_1^{\rm nw})&=\mathbf{x}, \qquad &\rho_1(T_1^{\rm ne})&=\mathbf{x}^{-1}, \qquad &\rho_1(T_1^{\rm sw})&=\mathbf{y},
\qquad   &\rho_1(T_1^{\rm ne})&=\mathbf{y}^{-1},  \\
\rho_2(T_2^{\rm ne})&=\mathbf{y}, \qquad &\rho_2(T_2^{\rm nw})&=\mathbf{y}^{-1}, \qquad &\rho_2(T_2^{\rm se})&=\mathbf{z},
\qquad   &\rho_2(T_2^{\rm sw})&=\mathbf{z}^{-1},  \\
&\ \qquad &\rho_3(T_3^{\rm nw})&=\mathbf{z}, \qquad &\rho_3(T_3^{\rm sw})&=\mathbf{x}^{-1}. &\ &\
\end{alignat*}
Define $\rho_4\in\mathcal{R}(T_4)$ by setting $\rho_4(\mathsf{a}')=\rho_3(\mathsf{a})$ for each directed arc $\mathsf{a}'$ of $T_4$. Then $\rho_4(\mathsf{c}')=\rho_4(\mathsf{a}')\lrcorner\rho_4(\mathsf{b}')$ for any directed arcs $\mathsf{a}',\mathsf{b}',\mathsf{c}'$ making a crossing in $T_4$ (see Figure \ref{fig:symmetric}), as $\rho_3(\mathsf{c})=\rho_3(\mathsf{a})\lrcorner\rho_3(\mathsf{b})$.
So indeed $\rho_4\in\mathcal{R}(T_4)$. It satisfies $\rho_4(T_4^{\rm ne})=\mathbf{z}^{-1}$, $\rho_4(T_4^{\rm se})=\mathbf{x}$.
We can glue the $\rho_i$'s into some $\rho\in\mathcal{R}_t^{\rm irr}(K)$ with $\rho(x_i)=\mathbf{x}$, $\rho(y_i)=\mathbf{y}$, $\rho(z_i)=\mathbf{z}$ for $i=1,2$.

The conjugacy classes of $(\mathbf{x},\mathbf{y},\mathbf{z})$ are parameterized by $b:={\rm tr}(\mathbf{x}\mathbf{z})$
and $r:={\rm tr}(\mathbf{x}\mathbf{y}\mathbf{z})$, subject to $f_t(a_1,a_2,b;r)=0$.
Therefore, $K$ admits a HDC.


\section{Examples}

In this section, we study three hyperbolic knots, explicitly describing the HDCs they admit.
Then we positively answer Question \ref{que:CD-2018}.

We clarify some basic facts first.

For $r=\mu+\mu^{-1}$, let
$$\omega_k(r)=\begin{cases} (\mu^k-\mu^{-k})/(\mu-\mu^{-1}),&\mu\ne\pm1 \\ k\mu^{k-1},&\mu=\pm1 \end{cases}.$$
By Cayley-Hamilton Theorem, if $\mathbf{x}\in G(r)$, then for each $k\in\mathbb{Z}$,  $$\mathbf{x}^k=\omega_k(r)\mathbf{x}-\omega_{k-1}(r)\mathbf{e}.$$
In particular, if $\mathbf{x}\in G(0)$, then $\mathbf{x}^2=-\mathbf{e}$; if $\mathbf{x}\in G(1)$, then $\mathbf{x}^3=-\mathbf{e}$.

\begin{lem}\label{lem:technique}
If $\mathbf{a},\mathbf{b}\in G(t)$ and ${\rm tr}(\mathbf{a}\mathbf{b})=1$, then
$\mathbf{a}\mathbf{b}\mathbf{a}^k\mathbf{b}\mathbf{a}=-\mathbf{b}^{k-2}$ for each $k\in\mathbb{Z}$;
in particular, $\mathbf{a}\mathbf{b}^2\mathbf{a}=-\mathbf{b}^{-2}$.
\end{lem}

\begin{proof}
Since ${\rm tr}(\mathbf{a}\mathbf{b}^2)={\rm tr}(\mathbf{a}(t\mathbf{b}-\mathbf{e}))=t\cdot{\rm tr}(\mathbf{a}\mathbf{b})-{\rm tr}(\mathbf{a})=0$,
we have $(\mathbf{a}\mathbf{b}^2)^2=-\mathbf{e}$, so that $\mathbf{a}\mathbf{b}^2\mathbf{a}=-\mathbf{b}^{-2}$.

For the general case, note that $(\mathbf{a}\mathbf{b})^3=-\mathbf{e}$, and then
\begin{align*}
\mathbf{a}\mathbf{b}\mathbf{a}^k\mathbf{b}\mathbf{a}\mathbf{b}^2
&=\mathbf{a}\mathbf{b}(\omega_k(t)\mathbf{a}-\omega_{k-1}(t)\mathbf{e})\mathbf{b}\mathbf{a}\mathbf{b}^2  \\
&=\omega_k(t)(\mathbf{a}\mathbf{b})^3\cdot\mathbf{b}-\omega_{k-1}(t)\mathbf{a}\mathbf{b}^2\mathbf{a}\cdot\mathbf{b}^2  \\
&=-\omega_k(t)\mathbf{b}+\omega_{k-1}(t)\mathbf{e}=-\mathbf{b}^k,
\end{align*}
implying $\mathbf{a}\mathbf{b}\mathbf{a}^k\mathbf{b}\mathbf{a}=-\mathbf{b}^{k-2}$.
\end{proof}

\begin{figure}[h]
  \centering
  \includegraphics[width=4.2cm]{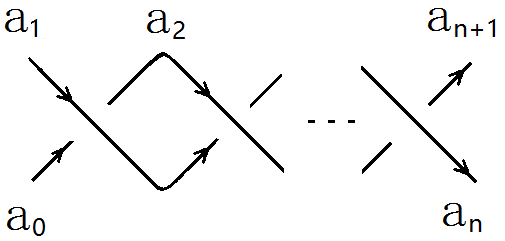}\\
  \caption{}\label{fig:tangle-rep}
\end{figure}

Consider the tangle $[n]$ with $n>0$, as shown in Figure \ref{fig:tangle-rep}.

Suppose $\tau\in\mathcal{R}_t([n])$, with $\tau(a_i)=\mathbf{a}_i$, and ${\rm tr}(\mathbf{a}_0\mathbf{a}_1)=s$.
Then $\mathbf{a}_{k+1}\mathbf{a}_k=\mathbf{a}_1\mathbf{a}_0$ for $1\le k\le n$, hence 
$\mathbf{a}_{k+2}=(\mathbf{a}_{k+1}\mathbf{a}_k)\lrcorner \mathbf{a}_k=(\mathbf{a}_1\mathbf{a}_0)\lrcorner\mathbf{a}_k$.
Consequently,
$$\mathbf{a}_3=\mathbf{a}_1\mathbf{a}_0\mathbf{a}_1(\mathbf{a}_1\mathbf{a}_0)^{-1}, \qquad
\mathbf{a}_4=(\mathbf{a}_1\mathbf{a}_0)^2\lrcorner\mathbf{a}_0.$$
As special cases, if $s=1$, then $(\mathbf{a}_1\mathbf{a}_0)^3=-\mathbf{e}$, so
$$\mathbf{a}_3=-\mathbf{a}_1\mathbf{a}_0\mathbf{a}_1(\mathbf{a}_1\mathbf{a}_0)^2
=-\mathbf{a}_1\mathbf{a}_0\mathbf{a}_1^2\mathbf{a}_0\mathbf{a}_1\cdot\mathbf{a}_0=\mathbf{a}_0;$$
if $s=0$, then $(\mathbf{a}_1\mathbf{a}_0)^2=-\mathbf{e}$, hence $\mathbf{a}_4=\mathbf{a}_0$.

By the formulas in \cite[Example 3.4]{Ch22}, with $h=1$, $r=t^2-s$,
\begin{align}
{\rm tr}(\mathbf{a}_1^{-1}\mathbf{a}_3)&=2+\grave{z}=2+(s+2-t^2)(s-2), \label{eq:formula-1} \\
{\rm tr}(\mathbf{a}_0^{-1}\mathbf{a}_3)&=2+z=2+(t^2-s-2)(s-1)^2.  \label{eq:formula-2}
\end{align}

\subsection{The $(3,3,3,4)$-pretzel knot}\label{sec:pretzel}

\begin{figure}[h]
  \centering
  \includegraphics[width=10cm]{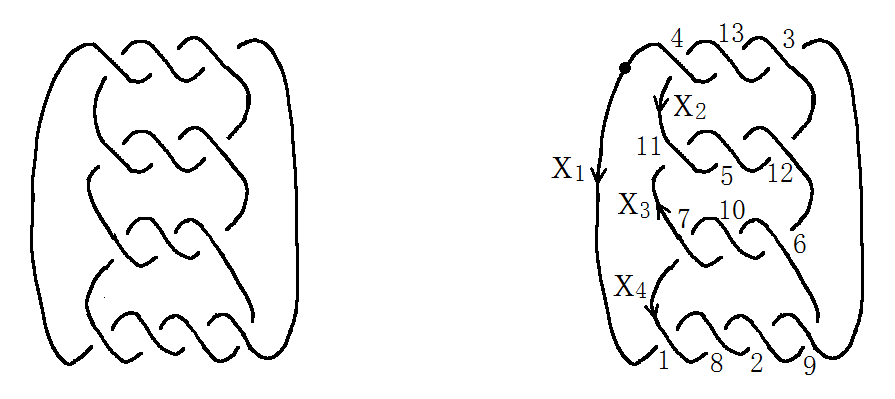}\\
  \caption{The $(3,3,3,4)$-pretzel knot.}\label{fig:pretzel}
\end{figure}

Let $P$ denote the $(3,3,3,4)$-pretzel knot, as shown in Figure \ref{fig:pretzel}.

Each quadruple $\underline{\mathbf{x}}=(\mathbf{x}_1,\mathbf{x}_2,\mathbf{x}_3,\mathbf{x}_4)\in G(t)^4$ satisfying
\begin{align}
{\rm tr}(\mathbf{x}_1\mathbf{x}_2)={\rm tr}(\mathbf{x}_2\mathbf{x}_3)={\rm tr}(\mathbf{x}_3\mathbf{x}_4)=1, \qquad
{\rm tr}(\mathbf{x}_4\mathbf{x}_1)=0 \label{eq:condition-0}
\end{align}
gives rise to a unique $\rho_{\underline{\mathbf{x}}}\in\mathcal{R}_t^{\rm irr}(K)$ with
$\mathbf{x}_i=\rho_{\underline{\mathbf{x}}}(x_i)$, $1\le i\le 4$ and $\mathbf{g}_{\rho_{\underline{\mathbf{x}}}}=\mathbf{e}$.

\begin{cla}\label{cla:HDC}
Let $\mathcal{Y}_t$ denote the set of conjugacy classes of $\underline{\mathbf{x}}\in G(t)^4$ satisfying {\rm(\ref{eq:condition-0})}. The map $\Psi_t:\mathcal{Y}_t\to\mathbb{C}^3$ defined by
$$\underline{\mathbf{x}}\mapsto({\rm tr}(\mathbf{x}_1\mathbf{x}_3),(\mathbf{x}_1\mathbf{x}_2\mathbf{x}_3),
{\rm tr}(\mathbf{x}_1\mathbf{x}_3\mathbf{x}_4))$$
is injective, and its image consists of $(t_{13},t_{123},t_{134})$ with $f_t(1,1,t_{13};t_{123})=f_t(1,0,t_{13};t_{134})=0$,
$t_{13}\ne t^2-2$, and $t_{13}\ne 2$ if $t^2\in\{1,2\}$.
\end{cla}

\begin{proof}
Given any $t_{13}\in\mathbb{B}_t$ and $t_{123}$ with $f_t(1,1,t_{13};t_{123})=0$, up to conjugacy there exists a unique
$(\mathbf{x}_1,\mathbf{x}_2,\mathbf{x}_3)\in G(t)^3$ such that ${\rm tr}(\mathbf{x}_1\mathbf{x}_2)={\rm tr}(\mathbf{x}_2\mathbf{x}_3)=1$, ${\rm tr}(\mathbf{x}_1\mathbf{x}_3)=t_{13}$, ${\rm tr}(\mathbf{x}_1\mathbf{x}_2\mathbf{x}_3)=t_{123}$. Since $t_{13}\in\mathbb{B}_t$, by Lemma \ref{lem:matrix} (ii-2), any $t_{134}$ with $f_t(1,0,t_{13};t_{134})=0$ further determines a unique $\mathbf{x}_4$ such that ${\rm tr}(\mathbf{x}_1\mathbf{x}_4)=0$,
${\rm tr}(\mathbf{x}_3\mathbf{x}_4)=1$, ${\rm tr}(\mathbf{x}_1\mathbf{x}_3\mathbf{x}_4)=t_{134}$.

From now on, let $t_{13}\in\{2,t^2-2\}$. Given any $t_{134}$ with $f_t(1,0,t_{13};t_{134})=0$, up to conjugacy there exits a
unique $(\mathbf{x}_1,\mathbf{x}_3,\mathbf{x}_4)\in G(t)^3$ such that
$${\rm tr}(\mathbf{x}_1\mathbf{x}_4)=0, \qquad {\rm tr}(\mathbf{x}_3\mathbf{x}_4)=1, \qquad {\rm tr}(\mathbf{x}_1\mathbf{x}_3\mathbf{x}_4)=t_{134}.$$ Since
${\rm tr}(\mathbf{x}_1\mathbf{x}_4)\ne {\rm tr}(\mathbf{x}_3\mathbf{x}_4)$, we have $\mathbf{x}_1\ne \mathbf{x}_3$.
\begin{itemize}
  \item If $\mathbf{x}_1=\mathbf{x}_3^{-1}$, then $t_{13}=2$, and $\mathbf{x}_1+\mathbf{x}_3=t\mathbf{e}$, so 
        $$1={\rm tr}(\mathbf{x}_1\mathbf{x}_4)+{\rm tr}(\mathbf{x}_3\mathbf{x}_4)={\rm tr}((\mathbf{x}_1+\mathbf{x}_3)\mathbf{x}_4)
        =t\cdot{\rm tr}(\mathbf{x}_4)=t^2.$$
  \item Conversely, if $(t^2,t_{13})=(1,2)$, then $f_t(1,0,t_{13};t_{134})=0$ forces $t_{134}=t$; since
        $(\mathbf{x}_3,\mathbf{x}_4)$ is irreducible and now ${\rm tr}(\mathbf{x}_1\mathbf{c})={\rm tr}(\mathbf{x}_3^{-1}\mathbf{c})$ for $\mathbf{c}\in\{\mathbf{x}_3,\mathbf{x}_4,\mathbf{x}_3\mathbf{x}_4\}$, by the uniqueness in Lemma \ref{lem:matrix} (ii-2), $\mathbf{x}_1=\mathbf{x}_3^{-1}$.
\end{itemize}
Hence $(t^2,t_{13})=(1,2)$ if and only if $\mathbf{x}_1=\mathbf{x}_3^{-1}$. When $(t^2,t_{13})=(1,2)$, there does not exist $\mathbf{x}_2$ with
${\rm tr}(\mathbf{x}_1\mathbf{x}_2)={\rm tr}(\mathbf{x}_2\mathbf{x}_3)=1$: otherwise, we would deduce
${\rm tr}(\mathbf{x}_1\mathbf{x}_2)+{\rm tr}(\mathbf{x}_2\mathbf{x}_3)=2\ne t^2={\rm tr}((\mathbf{x}_1+\mathbf{x}_3)\mathbf{x}_2)$.

Now suppose $(t^2,t_{13})\ne(1,2)$. Then $(\mathbf{x}_1,\mathbf{x}_3)$ is reducible and $\mathbf{x}_1\mathbf{x}_3\ne\mathbf{x}_3\mathbf{x}_1$.
Following the proof of Lemma \ref{lem:moduli}, in the context of Point 2 there, we can deduce $\varepsilon=-1$, so that $t_{13}=2$, and $t^2\ne 2$.
When these hold and given $t_{123}$ with $f_t(1,1,2;t_{123})=0$, there exists a unique $\mathbf{x}_2$ such that
${\rm tr}(\mathbf{x}_1\mathbf{x}_2)={\rm tr}(\mathbf{x}_2\mathbf{x}_3)=1$ and ${\rm tr}(\mathbf{x}_1\mathbf{x}_2\mathbf{x}_3)=t_{123}$.
\end{proof}

In conclusion, $P$ admits a HDC isomorphic to
\begin{align*}
\big\{(t,t_{13},t_{123},t_{134})\colon &t_{13}\ne t^2-2 \ \text{and\ }(t^2,t_{13})\ne(2,2),(1,2), \\
&t_{123}^2+t(t^2-2-t_{13})t_{123}+t_{13}^2+(1-t^2)t_{13}+t^2-2=0, \\
&t_{134}^2+t(t^2-1-t_{13})t_{134}+t_{13}^2-t^2t_{13}+2t^2-3=0\big\}.
\end{align*}

\begin{rmk}
\rm In \cite[Section 3.2]{Ch25}, we showed the existence of a HDC for a general 4-strand Montesinos knot.
Here for $P$, we precisely describe a HDC.
\end{rmk}

In general, given a knot $K$, each directed arc provides a meridional generator for $\pi(K)$.
Choose a preferred directed arc, which we always denote by $x_1$; let $\mathfrak{m}$ denote the corresponding meridian. Then the longitude $\mathfrak{l}$ paired with $\mathfrak{m}$ can be computed as follows.
Starting at the point marked with a bullet, walk along the diagram guided by the direction of $x_1$; denote the generator contributed by the over-crossing arc at the $k$-th crossing (labeled by $k$) by $\mathfrak{z}_k$. Then $\mathfrak{l}=\mathfrak{m}^{w(K)}\mathfrak{z}_1\cdots\mathfrak{z}_n$, where $n$ is the number of crossings, and $w(K)$ is the writhe whose presence ensures $\mathfrak{l}$ to represent $0$ in $H_1(E_K)$.

For the present knot, $w(P)=13$.

For $\underline{\mathbf{x}}\in \mathcal{Y}_t$, let $\mathbf{l}=\rho_{\underline{\mathbf{x}}}(\mathfrak{l})$.
As illustrated in the right part of Figure \ref{fig:pretzel}, $\rho_{\underline{\mathbf{x}}}(\mathfrak{z}_1)=\mathbf{x}_4^{-1}$, $\rho_{\underline{\mathbf{x}}}(\mathfrak{z}_2)=\mathbf{x}_1^{-1}\mathbf{x}_4^{-1}\mathbf{x}_1$, $\rho_{\underline{\mathbf{x}}}(\mathfrak{z}_3)=\mathbf{x}_2^{-1}$, etc. Hence
\begin{align*}
\mathbf{l}&=\mathbf{x}_1^{13}\cdot\mathbf{x}_4^{-1}\cdot\mathbf{x}_1^{-1}\mathbf{x}_4^{-1}\mathbf{x}_1\cdot\mathbf{x}_2^{-1}\cdot\mathbf{x}_1^{-1}
\cdot\mathbf{x}_2\mathbf{x}_3^{-1}\mathbf{x}_2^{-1}\cdot\mathbf{x}_4^{-1}\cdot\mathbf{x}_3^{-1} \\
&\ \ \ \ \cdot\mathbf{x}_4\mathbf{x}_1^{-1}\mathbf{x}_4^{-1}\cdot\mathbf{x}_1^{-1}\cdot\mathbf{x}_4\mathbf{x}_3^{-1}\mathbf{x}_4^{-1}
\cdot\mathbf{x}_2^{-1}\cdot\mathbf{x}_3^{-1}\cdot\mathbf{x}_2\mathbf{x}_1^{-1}\mathbf{x}_2^{-1}   \\
&=\mathbf{x}_1^{14}(\mathbf{x}_4\mathbf{x}_1)^{-2}\mathbf{x}_1\mathbf{x}_2^{-1}\mathbf{x}_1^{-1}\mathbf{x}_2\mathbf{x}_3^{-1}\mathbf{x}_2^{-1}
\mathbf{x}_4^{-1}\mathbf{x}_3^{-1}\mathbf{x}_4(\mathbf{x}_4\mathbf{x}_1)^{-2} \\
&\ \ \ \ \cdot\mathbf{x}_4^2\mathbf{x}_3^{-1}\mathbf{x}_4^{-1}\mathbf{x}_2^{-1}
\mathbf{x}_3^{-1}\mathbf{x}_2\mathbf{x}_1^{-1}\mathbf{x}_2^{-1}   \\
&=\mathbf{x}_1^{15}\mathbf{x}_2^{-1}\mathbf{x}_1^{-1}\mathbf{x}_2\mathbf{x}_3^{-1}\mathbf{x}_2^{-1}
(\mathbf{x}_4^{-1}\mathbf{x}_3^{-1}\mathbf{x}_4^3\mathbf{x}_3^{-1}\mathbf{x}_4^{-1})\mathbf{x}_2^{-1}
\mathbf{x}_3^{-1}\mathbf{x}_2\mathbf{x}_1^{-1}\mathbf{x}_2^{-1}   \\
&=-\mathbf{x}_1^{15}\mathbf{x}_2^{-1}\mathbf{x}_1^{-1}\mathbf{x}_2(\mathbf{x}_3^{-1}\mathbf{x}_2^{-1}
\mathbf{x}_3^5\mathbf{x}_2^{-1}\mathbf{x}_3^{-1})\mathbf{x}_2\mathbf{x}_1^{-1}\mathbf{x}_2^{-1}   \\
&=\mathbf{x}_1^{15}(\mathbf{x}_2^{-1}\mathbf{x}_1^{-1}\mathbf{x}_2^9\mathbf{x}_1^{-1}\mathbf{x}_2^{-1})   \\
&=-\mathbf{x}_1^{26}.
\end{align*}
In the third line, $(\mathbf{x}_4\mathbf{x}_1)^2=-\mathbf{e}$ is used, and in the last three lines, Lemma \ref{lem:technique} is applied.

Let $P(a/b)$ denote the Dehn filling of $E_P$ with slope $a/b$. For $\rho$ to descend to $\pi_1(P(a/b))\to G$, a sufficient and necessary condition is
$$\mathbf{e}=\mathbf{x}_1^a\mathbf{l}^b=(-1)^b\mathbf{x}_1^{a+26b},$$
which is equivalent to $t=\kappa+\kappa^{-1}$, with $\kappa^{a+26b}=(-1)^b$ and $\kappa\ne\pm1$.

By \cite[Theorem 3.6]{Wu96}, $P(a/b)$ is hyperbolic whenever $a/b\ne 1/0$. Thus, as long as $a+26b\ne\pm 1$, there exists $t\ne\pm2$ such that $\mathcal{X}^{\rm irr}(P(a/b))$ contains a 1-dimensional component contributed by $\mathcal{X}^{\rm irr}_t(P)$, so $P(a/b)$ provides a positive answer to Question \ref{que:CD-2018}. In particular, $P(0/1)$ does.

\subsection{Two knots with $12$ crossings}

Here we consider another two knots $Q_1,Q_2$, to illustrate the effects of the two conditional structures given in Section \ref{sec:conditional}.

Each knot also provides a positive answer to Question \ref{que:CD-2018}, for similar reasons as the last paragraphs in Section \ref{sec:pretzel}, which we no longer explain.

\begin{figure}[h]
  \centering
  \includegraphics[width=9.5cm]{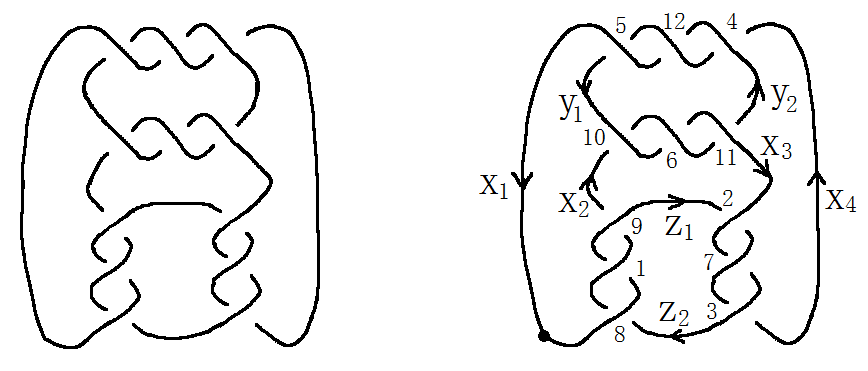}\\
  \caption{The knot $Q_1$.}\label{fig:Q1}
\end{figure}

Let $Q_1=D([3]\ast[3]\ast([-1/3]+[-1/3]))$, as shown in Figure \ref{fig:Q1}, left.

Each triple $\underline{\mathbf{w}}=(\mathbf{x},\mathbf{y},\mathbf{z})\in G(t)^3$ satisfying
\begin{align}
{\rm tr}(\mathbf{x}\mathbf{y})={\rm tr}(\mathbf{x}\mathbf{z})=1   \label{eq:condition-1}
\end{align}
determines a representation $\rho_{\underline{\mathbf{w}}}\in\mathcal{R}_t^{\rm irr}(Q_1)$ such that
$\rho_{\underline{\mathbf{w}}}(x_i)=\mathbf{x}$ for $1\le i\le 4$, and $\rho_{\underline{\mathbf{w}}}(y_i)=\mathbf{y}$, $\rho_{\underline{\mathbf{w}}}(z_i)=\mathbf{z}$ for $i=1,2$.

The conjugacy classes of $\underline{\mathbf{w}}\in G(t)^3$ with (\ref{eq:condition-1}) are parameterized by
$t_{23}={\rm tr}(\mathbf{y}\mathbf{z})$ and $t_{123}={\rm tr}(\mathbf{x}\mathbf{y}\mathbf{z})$, subject to $f_t(1,1,t_{23};t_{123})=0$.
Hence $Q_1$ admits a HDC isomorphic to
$$\big\{(t,t_{23},t_{123})\colon t_{123}^2+t(t^2-2-t_{23})t_{123}+t_{23}^2+(1-t^2)t_{23}+t^2-2=0\big\}.$$

Note that $w(Q_1)=12$.
Given $\underline{\mathbf{w}}$ with (\ref{eq:condition-1}), as illustrated in the right part of Figure \ref{fig:Q1}, we can compute
\begin{align*}
\rho_{\underline{\mathbf{w}}}(\mathfrak{l})&=\mathbf{x}^{12}\cdot \mathbf{x}\mathbf{z}^{-1}\mathbf{x}^{-1}\cdot \mathbf{x}^{-1}\cdot \mathbf{z}^{-1}\cdot \mathbf{y}^{-1}
\cdot \mathbf{x}^{-1}\cdot \mathbf{x}^{-1}\mathbf{y}^{-1}\mathbf{x}  \\
&\ \ \ \ \cdot \mathbf{x}\mathbf{z}^{-1}\mathbf{x}^{-1}\cdot \mathbf{x}^{-1}\cdot \mathbf{z}^{-1}\cdot \mathbf{y}^{-1}
\cdot\mathbf{x}^{-1}\cdot\mathbf{x}^{-1}\mathbf{y}^{-1}\mathbf{x}   \\
&=\mathbf{x}^{13}(\mathbf{z}^{-1}\mathbf{x}^{-2}\mathbf{z}^{-1})(\mathbf{y}^{-1}\mathbf{x}^{-2}\mathbf{y}^{-1})
\mathbf{x}^2(\mathbf{z}^{-1}\mathbf{x}^{-2}\mathbf{z}^{-1})(\mathbf{y}^{-1}\mathbf{x}^{-2}\mathbf{y}^{-1})\mathbf{x}  \\
&=\mathbf{x}^{13}(-\mathbf{x}^2)(-\mathbf{x}^2)\mathbf{x}^2(-\mathbf{x}^2)(-\mathbf{x}^2)\mathbf{x} \\
&=\mathbf{x}^{24}.
\end{align*}

\begin{figure}[H]
  \centering
  \includegraphics[width=9.5cm]{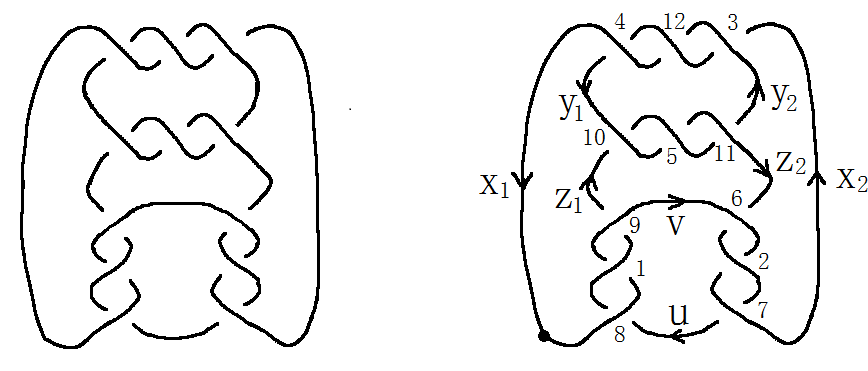}\\
  \caption{The knot $Q_2$.}\label{fig:Q2}
\end{figure}

Let $Q_2=D([3]\ast[3]\ast([-1/3]+[1/3]))$, as shown in Figure \ref{fig:Q2}, left.

Suppose $(\mathbf{x},\mathbf{y},\mathbf{z})\in G(t)^3$ satisfies
\begin{align}
{\rm tr}(\mathbf{x}\mathbf{y})={\rm tr}(\mathbf{y}\mathbf{z})=1, \qquad  {\rm tr}(\mathbf{y}\mathbf{z})\in\mathbb{B}_t.   \label{eq:condition-2-0}
\end{align}
Regard $[-1/3]$ as resulted from rotating $[3]$ by $-\pi/2$.
Given $s\in\mathbb{B}_t$, by (\ref{eq:formula-1}), (\ref{eq:formula-2}), necessary conditions for there to exist $\tau\in\mathcal{R}_t([-1/3])$ with $\tau(x_1)=\mathbf{x}$, $\tau(z_1)=\mathbf{z}$, $\tau(v)=\mathbf{v}$ and ${\rm tr}(\mathbf{z}\mathbf{v})=s$ are
\begin{align}
{\rm tr}(\mathbf{v}^{-1}\mathbf{x})&=2+(s+2-t^2)(s-2),  \label{eq:condition-2-1}  \\
{\rm tr}(\mathbf{z}^{-1}\mathbf{x})&=2-(s+2-t^2)(s-1)^2. \label{eq:condition-2-2}
\end{align}
Conversely, given $s\in\mathbb{B}_t$ satisfying (\ref{eq:condition-2-2}), we can take $\tilde{\mathbf{v}}\in G(t)$ with
${\rm tr}(\mathbf{z}\tilde{\mathbf{v}})=s$, and use $\mathbf{z},\tilde{\mathbf{v}}$ to construct a representation $\tilde{\tau}\in\mathcal{R}_t([-1/3])$ with $\tilde{\tau}(z_1)=\mathbf{z}$ and $\tilde{\tau}(v)=\tilde{\mathbf{v}}$; let $\tilde{\mathbf{x}}=\tilde{\tau}(x_1)$,
then by (\ref{eq:formula-1}), (\ref{eq:formula-2}),
\begin{align*}
{\rm tr}(\tilde{\mathbf{v}}^{-1}\tilde{\mathbf{x}})&=2+(s+2-t^2)(s-2),   \\
{\rm tr}(\mathbf{z}^{-1}\tilde{\mathbf{x}})&=2-(s+2-t^2)(s-1)^2={\rm tr}(\mathbf{z}^{-1}\mathbf{x}).
\end{align*}
The second equation ensures the existence of $\mathbf{c}\in G$ with $\mathbf{c}\lrcorner(\tilde{\mathbf{x}},\mathbf{z})=(\mathbf{x},\mathbf{z})$.
Then $\mathbf{v}:=\mathbf{c}\lrcorner\tilde{\mathbf{v}}$ fulfills ${\rm tr}(\mathbf{z}\mathbf{v})=s$ and (\ref{eq:condition-2-1}).


Each quadruple $\underline{\mathbf{w}}=(\mathbf{x},\mathbf{y},\mathbf{z},\mathbf{v})$ satisfying
${\rm tr}(\mathbf{z}\mathbf{v})=r$ and (\ref{eq:condition-2-0})--(\ref{eq:condition-2-2}) determines a representation $\rho_{\underline{\mathbf{w}}}\in\mathcal{R}_t^{\rm irr}(Q_2)$ with $\rho_{\underline{\mathbf{w}}}(x_i)=\mathbf{x}$, $\rho_{\underline{\mathbf{w}}}(y_i)=\mathbf{y}$, $\rho_{\underline{\mathbf{w}}}(z_i)=\mathbf{z}$ for $i=1,2$ and $\rho_{\underline{\mathbf{w}}}(v)=\mathbf{v}$.

The conjugacy classes of $\underline{\mathbf{w}}\in G(t)^4$ with (\ref{eq:condition-2-1}), (\ref{eq:condition-2-2}) are parameterized by
$t_{13}={\rm tr}(\mathbf{x}\mathbf{z})$ and $t_{123}={\rm tr}(\mathbf{x}\mathbf{y}\mathbf{z})$, subject to
$f_t(1,1,t_{13};t_{123})=0$. Hence $Q_2$ admits a HDC, which contains a Zariski open subset isomorphic to
\begin{align*}
\big\{(t,t_{13},t_{123},s)\colon &t_{123}^2+t(t^2-2-t_{13})t_{123}+t_{13}^2+(1-t^2)t_{13}+t^2-2=0, \\
&(s+2-t^2)(s-1)^2=t_{13}-t^2+2, \ t_{13},s\ne 2,t^2-2\big\}.
\end{align*}

Note that $w(Q_2)=6$.
For $\underline{\mathbf{w}}\in G(t)^4$ satisfying (\ref{eq:condition-2-1}), (\ref{eq:condition-2-2}), as illustrated in the right part of Figure \ref{fig:Q2}, we can compute
\begin{align*}
\rho_{\underline{\mathbf{w}}}(\mathfrak{l})&=\mathbf{x}^6\cdot \mathbf{x}\mathbf{u}^{-1}\mathbf{x}^{-1}\cdot \mathbf{x}\mathbf{u}\mathbf{x}^{-1}\cdot\mathbf{y}^{-1}
\cdot\mathbf{x}^{-1}\cdot\mathbf{y}\mathbf{z}^{-1}\mathbf{y}^{-1}\cdot \mathbf{v}  \\
&\ \ \ \ \cdot\mathbf{x}\cdot\mathbf{x}^{-1}\cdot\mathbf{v}^{-1}\cdot\mathbf{y}^{-1}\cdot\mathbf{z}^{-1}\cdot\mathbf{y}\mathbf{x}^{-1}\mathbf{y}^{-1}  \\
&=\mathbf{x}^6\mathbf{y}^{-1}\mathbf{x}^{-1}\mathbf{y}(\mathbf{z}^{-1}\mathbf{y}^{-2}\mathbf{z}^{-1})\mathbf{y}\mathbf{x}^{-1}\mathbf{y}^{-1}  \\
&=-\mathbf{x}^6(\mathbf{y}^{-1}\mathbf{x}^{-1}\mathbf{y}^4\mathbf{x}^{-1}\mathbf{y}^{-1}) \\
&=\mathbf{x}^{12}.
\end{align*}


\bigskip

\noindent
Haimiao Chen (orcid: 0000-0001-8194-1264)\ \ \ \ {\it chenhm@math.pku.edu.cn} \\
Department of Mathematics, Beijing Technology and Business University, \\
Liangxiang Higher Education Park, Fangshan District, Beijing, China.

\end{document}